\documentclass[11pt]{article}
\usepackage[margin=1in]{geometry}
\usepackage{amsmath,amssymb}
\usepackage{mathtools}
\usepackage{microtype}
\usepackage[hidelinks]{hyperref}
\AtBeginDocument{%
  \everydisplay\expandafter{\the\everydisplay
    \abovedisplayshortskip=\abovedisplayskip
    \belowdisplayshortskip=\belowdisplayskip
  }%
}
\title{The Quarter-Turn and the Ball}
\author{John Paul Steinberger}
\date{}

\begin{document}
\maketitle
\begin{abstract}
In his \emph{Tau Manifesto}, Michael Hartl presents a formula, credited
to Jeff Cornell, for the volume of the unit $n$-ball in terms of the
quarter-turn constant $\eta=\pi/2$. We read the formula as a set of clues
to a geometric proof, leading through orthant geometry and a
higher-dimensional version of Lambert's equal-area projection. We
conclude with some additional thoughts on the circle constant.
\end{abstract}

Let $\eta=\pi/2$ denote the radian measure of $90^\circ$. As
observed by Hartl in the \emph{Tau Manifesto}~\cite{hartl}, this constant gives
a particularly simple formula for the volume $V_n$ of the unit $n$-ball in terms
of elementary functions:
\begin{equation}\label{eq:whole}
 V_n=\frac{2^n\eta^{\lfloor n/2\rfloor}}{n!!}.
\end{equation}
For comparison, using $\tau=2\pi=4\eta$ gives
\begin{equation}\label{eq:tau}
 V_n=\frac{2^{\,n\bmod2}\tau^{\lfloor n/2\rfloor}}{n!!}.
\end{equation}
Here $n!!=n(n-2)(n-4)\cdots$, with $0!!=1!!=1$;
$\lfloor n/2\rfloor$ is $n/2$ rounded down; and $n\bmod2$ is $0$ or
$1$ according as $n$ is even or odd. 
(The similar formula with $\pi$ involves a factor of $2^{\lceil n/2\rceil}$ in the numerator.)

Hartl credits Jeff Cornell for \eqref{eq:whole}~\cite{hartl}.
While Hartl deduces \eqref{eq:whole}
from the classical formula for the volume of a ball, this note 
seeks 
to let the formula suggest its own proof.
We will
also compare the resulting proof with a short integral calculation.

The factor $2^n$ in \eqref{eq:whole} is the first clue: as Hartl notes, it counts 
the possible combinations of coordinate signs. Let $Q_n$ denote the
volume of the ball's intersection with the nonnegative orthant, where all coordinates are nonnegative.
Thus $Q_2$ is the area of a quarter-disk and $Q_3$ is the volume of
an eighth of a three-dimensional ball. Then \eqref{eq:whole} becomes
\begin{equation}\label{eq:formula}
 Q_n=\frac{\eta^{\lfloor n/2\rfloor}}{n!!}
\end{equation}
dividing by $2^n$.
The initial values are $Q_0=Q_1=1$: a point has zero-dimensional
volume $1$, while the nonnegative half of the unit $1$-ball has length
$1$. From there, \eqref{eq:formula} suggests proving the recurrence
\begin{equation}\label{eq:recurrence}
 Q_n=\frac{\eta}{n}Q_{n-2}\mathrlap{\qquad(n\geq2).}
\end{equation}

Next, the factor $1/n$ points to the conical relation between a ball's volume
and its surface area. Let $A_n$ denote the $(n-1)$-dimensional area of the
spherical part of the boundary of the orthant, excluding its flat
faces. Since the sphere has unit radius, $Q_n=A_n/n$. Thus
\begin{equation}\label{eq:area}
 A_n=\eta Q_{n-2}\mathrlap{\qquad(n\geq2)}
\end{equation}
is equivalent to \eqref{eq:recurrence}.

For $n=2$, \eqref{eq:area} amounts to the definition of $\eta$:
the unit quarter-circle has length $\eta$, and $Q_0=1$.
For $n=3$, we can verify \eqref{eq:area} manually using the familiar
equal-area property of Lambert's cylindrical projection~\cite{lambert}. This
projection sends each point of the unit sphere horizontally away
from the $z$-axis onto the surrounding cylinder of radius $1$, keeping
its height fixed.
The spherical octant maps to a quarter of the cylinder of height $1$,
whose area is $\eta\cdot1=\eta$.
This suggests that the proof of \eqref{eq:area} will rely on a higher-dimensional version of
Lambert's projection.

Split the coordinates of $\mathbb R^n$,
$n\geq3$, into a block of two and a block of $n-2$. Let
\[
 U=\{u\in\mathbb R^2:u_i\geq0,\ \|u\|=1\},\qquad
 W=\{v\in\mathbb R^{n-2}:v_i\geq0,\ \|v\|=1\}.
\]
Thus $U$ is the nonnegative quarter of the unit circle in $\mathbb R^2$,
whilst $W$ is the nonnegative orthant of a
unit sphere in the remaining coordinates. We also regard $u\in U$ and
$v\in W$ as vectors in the two perpendicular coordinate blocks of
$\mathbb R^n$. For $0\leq t\leq\eta$, define
\begin{gather*}
 \operatorname{Orth}(t)
   =(\cos(t)\,U)\times(\sin(t)\,W),\\
 \operatorname{Cyl}(t)
   =U\times(\sin(t)\,W).
\end{gather*}
We call these sets \emph{shells}. The shells $\operatorname{Orth}(t)$ partition
the spherical orthant, since $\cos(t)u+\sin(t)v$ is a vector of length
$1$ for each $u\in U$, $v\in W$, and every
point of the spherical orthant has this form for a unique $t$.
The shells $\operatorname{Cyl}(t)$ partition the product
of $U$ with the nonnegative orthant of the unit $(n-2)$-ball. We call this
product the cylinder. Since its $(n-1)$-dimensional area is $\eta Q_{n-2}$,
it remains to show that the cylinder and spherical
orthant have equal area.

For $0\leq t<\eta$, the correspondence
\[
 (\cos(t)u,\sin(t)v)\longleftrightarrow(u,\sin(t)v)
\]
maps $\operatorname{Orth}(t)$ bijectively onto $\operatorname{Cyl}(t)$.
It expands the first coordinate block to unit length
and leaves the second block fixed. For $n=3$, $W=\{1\}$, and we recover
Lambert's projection, with the shells as horizontal quarter-circles.

To show that this correspondence is area-preserving, we compare
thin slices between neighboring shells at $t$ and $t+dt$, accounting
for both shell area and slice width.
The width is measured between nearest points on the delimiting shells,
obtained by keeping $u\in U$, $v\in W$ fixed while $t$ varies.
These paths are perpendicular to the shells.

On the spherical orthant, fixing $u$ and $v$ and increasing $t$ to
$t+dt$ traces a unit-circle arc from $(\cos(t)u,\sin(t)v)$ to
$(\cos(t+dt)u,\sin(t+dt)v)$. Thus the slice
between $\operatorname{Orth}(t)$ and $\operatorname{Orth}(t+dt)$ has width
\[
 w_{\rm Orth}=dt
\]
since $u$ and $v$ are perpendicular vectors of length one.

On the cylinder, increasing $t$ to $t+dt$
traces a straight-line segment from $(u,\sin(t)v)$ to $(u,\sin(t+dt)v)$.
Thus the slice between $\operatorname{Cyl}(t)$ and
$\operatorname{Cyl}(t+dt)$ has width
\[
 w_{\rm Cyl}=\cos(t)\,dt
\]
to first order in $dt$, since $u$ stays fixed and the radius $\sin(t)$
increases at rate $\cos(t)$.

On the other hand, the shell-area ratio, cylinder to spherical
orthant, is $1/\cos(t)$, the expansion factor from
$\cos(t)U$ to $U$. This cancels the width ratio $\cos(t)$, so
corresponding slices have equal area to first order. Summing and passing
to the limit gives
\[
 A_n=\eta Q_{n-2},
\]
as required. This completes the reverse-engineered proof.

While the parity factor in \eqref{eq:tau} prevents reading that 
formula quite as simply,
the same proof does work for the whole ball: take $U$ to be
the full unit circle and $W$ the full unit sphere in the remaining
coordinates. The shell-area and width ratios are unchanged, while
the circle has length $\tau=4\eta$. Thus
\[
 V_n=\frac{\tau}{n}V_{n-2}.
\]
But the base cases now differ: $V_0=1$, whereas $V_1=2$,
accounting for the factor $2^{\,n\bmod2}$.

Conversely, other proofs for the whole ball can often be restricted to an
orthant. A familiar recursive integral gives a short example.
Use polar coordinates $(\rho,\theta)$ in the first two coordinates.
For fixed $\rho$, the remaining cross-section is a copy of the unit
$(n-2)$-ball orthant scaled by $\sqrt{1-\rho^2}$, while the volume
element in the first two coordinates is $\rho\,d\rho\,d\theta$. Thus
\[
 Q_n=Q_{n-2}\int_0^\eta d\theta
       \int_0^1(1-\rho^2)^{(n-2)/2}\rho\,d\rho
     =\frac{\eta}{n}Q_{n-2},
\]
where the last integral is evaluated by putting $s=1-\rho^2$.

In retrospect,
$\eta$'s advantage over $\tau$ in the ball-volume formula just 
comes down to the base cases: $Q_0=Q_1$, whereas $V_0\ne V_1$.
This seems almost anecdotical.
Still, at this point we would like to 
offer a sincere, if heretical, proposal: $\eta$ 
\emph{should} be the circle constant.
Basic trigonometric identities such as
\begin{gather*}
 \sin(t+\eta)=\cos(t),\\
 \cos(t+\eta)=-\sin(t),\\
 \tan(t+\eta)=-1/\tan(t)
\end{gather*}
would become friendlier, 
and the horizontal axes of trigonometry graphs could
shed many a fraction. 
More to the point, of
$90^\circ$, $180^\circ$, and $360^\circ$, the universe seems far more
interested in one than in the other two.

There is also the rational
approximation $\eta\approx11/7$, compared with $\pi\approx22/7$ and
$\tau\approx44/7$.
We have difficulty believing that $11/7$ was meant to be doubled this whole time.

For higher-dimensional sphere--cylinder area identities, see
Rudin~\cite{rudin} and Coll--Dodd--Harrison~\cite{coll-dodd-harrison}.
For other treatments of ball volumes, see Badger~\cite{badger}; for
generalized unit balls, see Wang~\cite{wang}.

\end{document}